# SYMMETRIC PAIRS OF UNBOUNDED OPERATORS IN HILBERT SPACE, AND THEIR APPLICATIONS IN MATHEMATICAL PHYSICS

PALLE E. T. JORGENSEN AND ERIN P. J. PEARSE

ABSTRACT. In a previous paper, the authors introduced the idea of a symmetric pair of operators as a way to compute self-adjoint extensions of symmetric operators. In brief, a symmetric pair consists of two densely defined linear operators $A$ and $B$, with $A \subseteq B^\star$ and $B \subseteq A^\star$. In this paper, we will show by example that symmetric pairs may be used to deduce closability of operators and sometimes even compute adjoints. In particular, we prove that the Malliavin derivative and Skorokhod integral of stochastic calculus are closable, and the closures are mutually adjoint. We also prove that the basic involutions of Tomita-Takesaki theory are closable and that their closures are mutually adjoint. Applications to functions of finite energy on infinite graphs are also discussed, wherein the Laplace operator and inclusion operator form a symmetric pair.

## 1. INTRODUCTION

In this paper, we discuss a number of questions in infinite-dimensional analysis: in stochastic calculus for Gaussian fields, and in von Neumann algebra theory. While these settings initially appear rather disparate, we show that they are both amenable to certain spectral theoretic considerations once one has identified an appropriate Hilbert space theoretic framework (i.e., appropriate unbounded linear operators with dense domain in Hilbert space). It turns out to be essential that these operators are *closable* (see Definition 2.4 below). Rather than studying just a single operator, we shall address these questions with the use of a carefully chosen *symmetric pair* of operators; see Definition 2.1, below.

The mechanism of a symmetric pair of operators $(A, B)$ may be used to provide short and elegant proofs of closability, and may sometimes allow for computation of the adjoint. The idea is as follows: given Hilbert spaces $\mathcal{H}_1$ and $\mathcal{H}_2$ and an unbounded linear operator $A$ with dense domain in $\mathcal{H}_1$, one can check the closability of $A : \mathcal{H}_1 \to \mathcal{H}_2$ by constructing an operator $B : \mathcal{H}_2 \to \mathcal{H}_1$ with the property that $A \subseteq B^\star$ and $B \subseteq A^\star$. Given such a pair, it has been shown in [JP10, Lem. 2.2] that both $A$ and $B$ are closable. Furthermore, in two of the applications presented below (stochastic calculus in §3, and Tomita-Takesaki theory in §4), if the operator $B$ in the symmetric pair is specified carefully, one may even be able to show that $\overline{A} = B^\star$ and $\overline{B} = A^\star$ (in which case we say the symmetric pair is maximal; cf. Definition 2.11). Therefore, this approach can be used to find the adjoint of an unbounded operator. We demonstrate the efficacy of this approach in several contexts. We also provide an example of a symmetric pair which is not maximal (the graph Laplacian in the energy space of an infinite network) in §5.

1.1. **Gaussian processes.** By *Gaussian fields*, we refer to a general family of Gaussian processes which are naturally indexed by a fixed Hilbert space; see Definition 3.1 for details. In this setting, there is an infinite-dimensional calculus of variations (often called Malliavin calculus), and the pair of operators we consider consists of a differential and an integral operator.

The differential operator is often called the Malliavin derivative, and the other is called the stochastic integral, or Skorokhod integral. We refer to [Hid, Bel, PS, Gaw] for details.

In this context, we prove that the Malliavin derivative and Skorokhod integral are closable, and prove that the closures of these operators are mutuallly adjoint. Closability of the Malliavin derivative is fundamental to the Malliavin calculus.

## 1.2. Tomita-Takesaki theory.

Tomita-Takesaki theory was devised in the 1970s as a tool for understanding von Neumann algebras which are type III factors; cf. [Tak1]. A *factor* is a von Neumann algebra whose center consists only of multiples of the identity operator. Given a factor $\mathfrak{M}$, a rough outline of the classification is as follows

- Type I. $B(\mathcal{H})$, the algebra of all bounded linear operators on some Hilbert space $\mathcal{H}$. It is called type $I_n$ if $\dim \mathcal{H} = n$.
- Type $II_1$. $\mathfrak{M}$ is of type $II_1$ iff $\mathfrak{M}$ is not type I, but has a finite and faithful trace.
- Type $II_\infty$. $\mathfrak{M}$ is of type $II_\infty$ iff $\mathfrak{M}$ is not type $II_1$, but has an infinite but semifinite faithful trace. Here, *semifinite* means that the domain is dense in $\mathfrak{M}$ in the $w^*$ topology.
- Type III. Everything else. So no trace of a trace.

The study of type III von Neumann algebras $\mathfrak{M}$ and associated commutant $\mathfrak{M}'$ goes back the early days. Among the big advances from Tomita-Takesaki theory is the following result. Given $\mathfrak{M}$ type III and with a choice of a faithful normal state. Then we can make explicit the assertion that the commutant $\mathfrak{M}'$ is realized as an anti-isomorphic copy of $\mathfrak{M}$, even spatially implemented. Here, *spatially* refers to an operator in $\mathcal{H}$. In this case, the operator is called the modular conjugation and denoted $J$, and so one has $J\mathfrak{M}J = \mathfrak{M}'$; cf. Remark 4.7 and Remark 4.9. We refer to [KR1, BR, Kad, Tak1, Tak2, vD, HV] for details; see [DS] for more general background.

In this context, we prove that the basic involutions of Tomita-Takesaki theory are closable and compute the adjoints of their closures (cf. [KR1, §9.2]). Closability of the involutions of Tomita-Takesaki theory are fundamental for constructing the modular automorphisms necessary for the study of von Neumann algebras which are type III factors.

While closability of the pertinent operators is previously known in both of these example scenarios, we found in each case that the method of symmetric pairs provides for shorter, simpler, and more elegant proofs. It yields closability and other conclusions in a simple and unified manner, and serves to unify two otherwise disparate areas of infinite-dimensional analysis. We invite the reader to compare the relatively brief and straightforward proofs of Theorem 4.6 and Theorem 4.11 with the machinery of [KR1, Thm. 2.7.8, Rem. 5.6.3, Lem. 9.2.1, Lem. 9.2.28, Cor. 9.2.30, Cor. 9.2.31].

## 2. Symmetric pairs

In a previous paper, the authors introduced the idea of a symmetric pair of operators as a way to compute self-adjoint extensions of symmetric operators; see [JP10, JPT]. In this paper, we will show by example that symmetric pairs may be used to deduce closability of operators and sometimes even compute adjoints. See [DS, BR] for background.

**Definition 2.1.** Suppose $\mathcal{H}_1$ and $\mathcal{H}_2$ are Hilbert spaces and $A, B$ are operators with dense domains $\dom A \subseteq \mathcal{H}_1$ and $\dom B \subseteq \mathcal{H}_2$ and

$$A : \dom A \subseteq \mathcal{H}_1 \to \mathcal{H}_2 \qquad \text{and} \qquad B : \dom B \subseteq \mathcal{H}_2 \to \mathcal{H}_1.$$

If $A$ and $B$ are linear operators, we say that $(A, B)$ is a *symmetric pair* iff

$$\langle A\varphi, \psi \rangle_{\mathcal{H}_2} = \langle \varphi, B\psi \rangle_{\mathcal{H}_1}, \qquad \text{for all } \varphi \in \dom A, \psi \in \dom B. \tag{2.1}$$

In other words, $(A, B)$ is a symmetric pair iff

$$A \subseteq B^\star \qquad \text{and} \qquad B \subseteq A^\star. \tag{2.2}$$

In the case when $A$ and $B$ are conjugate linear operators (see §4), we say that $(A, B)$ is a *symmetric pair* iff

$$\langle A\varphi, \psi \rangle_{\mathcal{H}_2} = \overline{\langle \varphi, B\psi \rangle}_{\mathcal{H}_1}, \qquad \text{for all } \varphi \in \operatorname{dom} A, \psi \in \operatorname{dom} B. \tag{2.3}$$

*Remark* 2.2 (Integration by parts). The best understood symmetric pairs arise from some form of integration by parts. For example, consider $\mathcal{H}_1 = \mathcal{H}_2 = L^2(\mathfrak{g})$, where $\mathfrak{g}$ is the standard Gaussian measure (i.e., $N(0,1)$) on $\mathbb{R}$: $d\mathfrak{g}(x) := (2\pi)^{-1/2} e^{x^2/2} dx$. In this case, one can easily check that $(A, B)$ is a symmetric pair for the following two operators:

$$\operatorname{dom} A := \{f \in L^2(\mathfrak{g}) : f'(x) \in L^2(\mathfrak{g})\}, \qquad A := \tfrac{d}{dx} \tag{2.4}$$

$$\operatorname{dom} B := \{f \in L^2(\mathfrak{g}) : xf(x) - f'(x) \in L^2(\mathfrak{g})\}, \qquad B := -\tfrac{d}{dx} + M_x, \tag{2.5}$$

where $M_x$ denotes the operation of multiplying by the independent variable $x$. Both domains are clearly dense in $L^2(\mathfrak{g})$. Below, we give examples of important symmetric pairs which are *not* related by integration by parts.

**Definition 2.3.** The *graph* of a linear operator $A: \mathcal{H}_1 \to \mathcal{H}_2$ is

$$\Gamma(A) := \left\{ \begin{bmatrix} x \\ y \end{bmatrix} \in \mathcal{H}_1 \oplus \mathcal{H}_2 : y = Ax \right\}.$$

**Definition 2.4.** A densely defined linear operator

$$A: \mathcal{H}_1 \to \mathcal{H}_2, \qquad \operatorname{dom} A \subseteq \mathcal{H}_1,$$

is said to be *closable* iff $\overline{\Gamma(A)}$ (the closure of $\Gamma(A)$ in the natural Hilbert norm of $\mathcal{H}_1 \oplus \mathcal{H}_2$) is the graph of an operator. In this case, the *closure* of $A$ is denoted $\overline{A}$, and one has $\overline{\Gamma(A)} = \Gamma(\overline{A})$.

**Theorem 2.5** ([JP10, Lem. 2.2]). *If $(A, B)$ is a symmetric pair, then $A$ and $B$ are each closable operators. Moreover,*

*(1) $A^\star \overline{A}$ is densely defined and self-adjoint with $\operatorname{dom} A^\star \overline{A} \subseteq \operatorname{dom} \overline{A} \subseteq \mathcal{H}_1$, and*
*(2) $B^\star \overline{B}$ is densely defined and self-adjoint with $\operatorname{dom} B^\star \overline{B} \subseteq \operatorname{dom} \overline{B} \subseteq \mathcal{H}_2$.*

*Remark* 2.6. Note that the proof of Theorem 2.5 given in [JP10, Lem. 2.2] is only for the linear case (2.1). However, the same proof works, mutatis mutandis, for the conjugate linear case (2.3).

**Corollary 2.7.** *If $(A, B)$ is a symmetric pair, then there is a partial isometry $V: \mathcal{H}_1 \to \mathcal{H}_2$ such that $\overline{A} = V(A^\star \overline{A})^{1/2} = (B^\star \overline{B})^{1/2} V$. In particular,*

$$\operatorname{spec}_{\mathcal{H}_1}(A^\star \overline{A}) \setminus \{0\} = \operatorname{spec}_{\mathcal{H}_2}(B^\star \overline{B}) \setminus \{0\}.$$

*Let $P_1$ denote the orthogonal projection onto $\mathcal{H}_1 \ominus \ker A$, and let $P_2$ denote the projection onto $\overline{\operatorname{ran} A}$ in $\mathcal{H}_2$. Then $V^\star V = P_1$ and $VV^\star = P_2$.*

Corollary 2.7 is an immediate consequence of Theorem 2.5.

*Remark* 2.8. The properties of the partial isometry $V$ in Corollary 2.7 are an important consequence of Theorem 2.5. In the two applications of symmetric pairs given in this paper (the Malliavin derivative §3 and the modular conjugation of Tomita-Takesaki theory in §4), we make use of $V$, and we are able to compute it.

- We show that in the case of the Malliavin derivative (cf. Corollary 3.18), the defect of the partial isometry is one-dimensional. More precisely, the partial isometry $V$ is isometric except for one dimension in the case of the Malliavin derivative, and the one-dimensional space is spanned by the constant function **1** on the sample space $\Omega$.
- We show that in the case of modular theory, the defect of the partial isometry is zero (cf. Theorem 4.10), so that $V$ is, in fact, an isometry, called the modular conjugation.

All of these results follow from our use of the particular symmetric pair at hand. See [DS, KR2] for background.

*Remark* 2.9. It follows from Theorem 2.5 that (2.2) may be rewritten as

$$\overline{A} \subseteq B^\star \qquad \text{and} \qquad \overline{B} \subseteq A^\star, \tag{2.6}$$

which will be useful for Lemma 2.15.

**Lemma 2.10.** *Let $(A, B)$ be a symmetric pair. Then $\overline{A} = B^\star$ iff $\overline{B} = A^\star$.*

*Proof.* Due to the symmetry of the situation, it suffices to show only one direction, so assume $\overline{A} = B^\star$. Upon taking adjoints, we have $A^\star = (\overline{A})^\star = B^{\star\star} = \overline{B}$, and the proof is complete. □

**Definition 2.11.** A symmetric pair $(A, B)$ is called *maximal* iff equality holds upon taking closures in (2.2), i.e., iff the closure of one operator is the adjoint of the other:

$$\overline{A} = B^\star \qquad \text{and} \qquad \overline{B} = A^\star. \tag{2.7}$$

By Lemma 2.10, only one of the equalities in (2.7) needs to be verified in a given instance.

*Remark* 2.12 (Examples of maximal symmetric pairs). The Malliavin derivative in §3 and the polar decomposition of Tomita-Takesaki theory in §4 below will furnish examples of maximal symmetric pairs. For an example of a symmetric pair which is not maximal; see §5, especially Remark 5.19.

**Definition 2.13.** In Corollary 2.14 and Corollary 2.7 and elsewhere, we use the notation $\mathcal{H} \ominus \mathcal{K}$ to denote the orthogonal complement of $\mathcal{K}$ in $\mathcal{H}$:

$$\mathcal{H} \ominus \mathcal{K} := \{h \in \mathcal{H} \,\vdots\, \langle h, k \rangle_\mathcal{H} = 0, \text{ for all } k \in \mathcal{K}\}. \tag{2.8}$$

**Corollary 2.14.** *For any symmetric pair $(A, B)$,*

$$\Gamma(B^\star) \ominus \Gamma(A) \neq \{0\} \text{ in } \mathcal{H}_1 \oplus \mathcal{H}_2 \quad \Longleftrightarrow \quad \Gamma(A^\star) \ominus \Gamma(B) \neq \{0\} \text{ in } \mathcal{H}_1 \oplus \mathcal{H}_2.$$

**Lemma 2.15.** *Suppose $(A, B)$ is a symmetric pair. The containment $\overline{A} \subseteq B^\star$ of (2.6) is strict if and only if there is a nonzero $\psi$ satisfying*

$$\psi \in \text{dom}(A^\star B^\star) \qquad \text{and} \qquad A^\star B^\star \psi = -\psi. \tag{2.9}$$

*Proof.* Recall that $A^\star = (\overline{A})^\star$ and that $A^{\star\star} = \overline{A}$. The containment $A \subseteq B^\star$ is strict if and only if there is a nonzero $\psi \in \text{dom } B^\star$ such that

$$0 = \left\langle \begin{bmatrix} u \\ Au \end{bmatrix}, \begin{bmatrix} \psi \\ B^\star \psi \end{bmatrix} \right\rangle = \langle u, \psi \rangle_{\mathcal{H}_1} + \langle Au, B^\star \psi \rangle_{\mathcal{H}_2}, \qquad \text{for all } u \in \text{dom } A, \tag{2.10}$$

where the first equality follows from the definition of containment of operators, and the second by the definition of the usual inner product on $\mathcal{H}_1 \oplus \mathcal{H}_2$. Observe that (2.10) implies $B^\star \psi \in \text{dom } A^\star$ and hence $\psi \in \text{dom}(A^\star B^\star)$, so that (2.10) is seen to be equivalent to (2.9). □

2.1. **The relation between symmetric operators and symmetric pairs.** In the case of a symmetric operator, the question of existence of proper symmetric extensions (or possibly even proper self-adjoint extensions) is decided by certain eigenspaces of the adjoint operator $A^\star$. Lemma 2.15 provides an analogous result for symmetric pairs, and we now extend this fruitful analogy.

The notion of "symmetric pairs of unbounded operators" generalizes the concept from von Neumann of a single symmetric operator as in Definition 2.16. In Theorem 2.17, we make the relation between symmetric operators and symmetric pairs more precise. Theorem 2.17 shows that for every symmetric pair of operators $(A, B)$, there is an associated canonical symmetric operator, denoted $L$ below, with $L$ acting in the direct sum Hilbert space, and determined uniquely by the given pair $(A, B)$. We further show in Corollary 2.20 that the notion of "deficiency" for a pair $(A, B)$ is directly connected to the deficiency of the symmetric operator $L$ (in the sense of von Neumann), and hence to von Neumann deficiency indices.

**Definition 2.16.** Let $\mathcal{K}$ be a Hilbert space and let $\mathcal{D} \subseteq \mathcal{K}$ be a dense subspace and consider a linear operator $L$ on $\mathcal{K}$ with $\text{dom } L = \mathcal{D}$. We say that $L$ is *symmetric* (or *Hermitian*) iff

$$\langle Lw, z \rangle_\mathcal{K} = \langle w, Lz \rangle_\mathcal{K}, \qquad \text{for all } w, z \in \mathcal{D}, \tag{2.11}$$

or equivalently, iff $L \subseteq L^\star$.

**Theorem 2.17.** *Suppose $\mathcal{H}_1$ and $\mathcal{H}_2$ are Hilbert spaces and $A, B$ are operators with dense domains $\operatorname{dom} A \subseteq \mathcal{H}_1$ and $\operatorname{dom} B \subseteq \mathcal{H}_2$ and*

$$A : \operatorname{dom} A \subseteq \mathcal{H}_1 \to \mathcal{H}_2 \quad \text{and} \quad B : \operatorname{dom} B \subseteq \mathcal{H}_2 \to \mathcal{H}_1.$$

*Define $\mathcal{K} := \mathcal{H}_1 \oplus \mathcal{H}_2$ and let $L : \mathcal{K} \to \mathcal{K}$ be the densely defined linear operator in $\mathcal{K}$ given by $L(u \oplus v) := Bv \oplus Au$ for $u \in \operatorname{dom} A$ and $v \in \operatorname{dom} B$. In other words,*

$$\operatorname{dom} L := \operatorname{dom} A \oplus \operatorname{dom} B, \quad \text{and} \quad L := \begin{bmatrix} 0 & B \\ A & 0 \end{bmatrix}. \tag{2.12}$$

*Then $(A, B)$ is a symmetric pair iff $L$ is symmetric.*

*Proof.* Let $w = u \oplus v \in \operatorname{dom} A \oplus \operatorname{dom} B$ and let $z = x \oplus y \in \operatorname{dom} A \oplus \operatorname{dom} B$. On the one hand,

$$\langle Lw, z \rangle_\mathcal{K} = \left\langle \begin{bmatrix} Bv \\ Au \end{bmatrix}, \begin{bmatrix} x \\ y \end{bmatrix} \right\rangle = \langle Bv, x \rangle_{\mathcal{H}_1} + \langle Au, y \rangle_{\mathcal{H}_2}. \tag{2.13}$$

On the other hand,

$$\langle w, Lz \rangle_\mathcal{K} = \left\langle \begin{bmatrix} u \\ v \end{bmatrix}, \begin{bmatrix} By \\ Ax \end{bmatrix} \right\rangle = \langle v, Ax \rangle_{\mathcal{H}_2} + \langle u, By \rangle_{\mathcal{H}_1}. \tag{2.14}$$

Now from (2.13) and (2.14), it is clear that $\langle Lw, z \rangle_\mathcal{K} = \langle w, Lz \rangle_\mathcal{K}$ iff

$$\langle Bv, x \rangle_{\mathcal{H}_1} + \langle Au, y \rangle_{\mathcal{H}_2} = \langle v, Ax \rangle_{\mathcal{H}_2} + \langle u, By \rangle_{\mathcal{H}_1},$$

which holds for all $x, u \in \operatorname{dom} A$ and $y, v \in \operatorname{dom} B$ precisely when (2.1) is true. $\square$

**Corollary 2.18.** *For $L$ defined as in (2.12), $L^\star(x \oplus y) = A^\star y \oplus B^\star x$ for $y \in \operatorname{dom} A^\star$ and $x \in \operatorname{dom} B^\star$. In other words,*

$$\operatorname{dom} L^\star = \operatorname{dom} B^\star \oplus \operatorname{dom} A^\star, \quad \text{and} \quad L^\star = \begin{bmatrix} 0 & A^\star \\ B^\star & 0 \end{bmatrix} \tag{2.15}$$

*Proof.* The proof of (2.15) follows by applying (2.12) to obtain

$$\left\langle \begin{bmatrix} x \\ y \end{bmatrix}, L \begin{bmatrix} u \\ v \end{bmatrix} \right\rangle = \left\langle \begin{bmatrix} A^\star y \\ B^\star x \end{bmatrix}, \begin{bmatrix} u \\ v \end{bmatrix} \right\rangle. \qquad \square$$

**Definition 2.19.** For a densely defined symmetric operator $T : \mathcal{K} \to \mathcal{K}$, the *defect spaces* are

$$\operatorname{def}_\pm(T) := \operatorname{Eig}_{\pm \mathbbm{i}}(T^\star) = \{u_\pm \in \operatorname{dom} T^\star : T^\star u_\pm = \pm \mathbbm{i} u_\pm\}. \tag{2.16}$$

Also, the *deficiency indices* of $T$ are

$$n_\pm(T) := \dim \operatorname{def}_\pm(T). \tag{2.17}$$

It is well-known, due to the work of von Neumann, that a densely defined symmetric operator has self-adjoint extensions if and only if its deficiency indices are equal ($n_+(T) = n_-(T)$) and that it has a *unique* self-adjoint extension (i.e., is *essentially self-adjoint*) iff both deficiency indices equal 0.

**Theorem 2.20.** *Let $(A, B)$ be a asymmetric pair and let $L$ be defined as in (2.12).*
*(a) The deficiency indices of $L$ are equal: $n_+(L) = n_-(L)$.*
*(b) Furthermore, $L$ is essentially self-adjoint iff $(A, B)$ is a maximal symmetric pair.*

*Proof.* (a) If $u \oplus v$ is an eigenvector of $L^\star$ with eigenvalue $\mathbbm{i}$, then (2.15) gives

$$\begin{bmatrix} A^\star v \\ B^\star u \end{bmatrix} = L^\star \begin{bmatrix} u \\ v \end{bmatrix} = \mathbbm{i} \begin{bmatrix} u \\ v \end{bmatrix} = \begin{bmatrix} \mathbbm{i} u \\ \mathbbm{i} v \end{bmatrix}, \tag{2.18}$$

which is equivalent to

$$v \in \operatorname{dom} A^\star, \text{ with } A^\star v = \mathbbm{i} u \quad \text{and} \quad u \in \operatorname{dom} B^\star, \text{ with } B^\star u = \mathbbm{i} v. \tag{2.19}$$

Now can see that

$$\begin{bmatrix} u \\ v \end{bmatrix} \in \operatorname{def}_+(L) \quad \Longleftrightarrow \quad \begin{bmatrix} -u \\ v \end{bmatrix} \in \operatorname{def}_-(L). \tag{2.20}$$

by the following computation:

$$L^\star \begin{bmatrix} -u \\ v \end{bmatrix} \overset{=}{_{(i)}} \begin{bmatrix} A^\star v \\ -B^\star u \end{bmatrix} \overset{=}{_{(ii)}} \begin{bmatrix} \mathfrak{i} u \\ -\mathfrak{i} v \end{bmatrix} = -\mathfrak{i} \begin{bmatrix} -u \\ v \end{bmatrix},$$

where $(i)$ follows from (2.15) and $(ii)$ follows from (2.19). The other computation is similar. This gives a bijection between $\mathrm{def}_+(L)$ and $\mathrm{def}_-(L)$, whence $n_+(L) = n_-(L)$.

(b) By the theory of von Neumann, a self-adjoint extension of $L$ corresponds to a partial isometry $Q$ which maps the defect space $\mathrm{def}_+(L)$ onto the defect space $\mathrm{def}_-(L)$. We denote this extension by $L_Q$ and observe that $L \subseteq L_Q \subseteq L^\star$, whence (2.12) and (2.15) give the representation

$$L_Q = \begin{bmatrix} 0 & B_Q \\ A_Q & 0 \end{bmatrix} \tag{2.21}$$

for some operators $A_Q$ and $B_Q$ satisfying $A \subseteq A_Q \subseteq B^\star$ and $B \subseteq B_Q \subseteq A^\star$. Taking adjoints yields $\overline{B} \subseteq A_Q^\star \subseteq A^\star$ and $\overline{A} \subseteq B_Q^\star \subseteq B^\star$, and since the self-adjointness of $L_Q$ implies $A_Q^\star = B_Q$ and $B_Q^\star = A_Q$, we have

$$\overline{A} \subseteq A_Q \subseteq B^\star \quad \text{and} \quad \overline{B} \subseteq B_Q \subseteq A^\star. \tag{2.22}$$

From (2.22), we see that maximality of the symmetric pair implies $L_Q = \overline{L}$ and conversely. $\square$

*Remark* 2.21 (Self-adjoint extensions of $L$). From Corollary 2.20 and its proof, it is clear that every self-adjoint extension of $L$ corresponds to an extension $(A_Q, B_Q)$ of $(A, B)$ in the sense that $(A_Q, B_Q)$ is a symmetric pair satisfying $A \subseteq A_Q$ and $B \subseteq B_Q$.

*Remark* 2.22. Combining (2.18) with (2.20), one can see that $u \oplus v \in \mathrm{def}_+(L^\star)$ if and only if

$$\begin{bmatrix} u \\ v \end{bmatrix} = \begin{bmatrix} u \\ -\mathfrak{i} B^\star u \end{bmatrix} = \begin{bmatrix} \mathfrak{i} A^\star v \\ v \end{bmatrix}, \tag{2.23}$$

where the appropriate domain assumptions for either direction of the implication. The corresponding identity holds for $\mathrm{def}_-(L^\star)$, and both will be useful in the sequel.

According to von Neumann's theory of self-adjoint extensions of a symmetric operator $T$, the set of self-adjoint extensions of $T$ is parametrized by partial isometries from $\mathrm{def}_+(T)$ to $\mathrm{def}_-(T)$. Denoting such a map by $Q$ (as in the proof of Theorem 2.20(b)) and the Cayley transform by $C(T) := (\mathfrak{i} - T)(\mathfrak{i} + T)^{-1}$, we apply this to $T = L$ and $\mathcal{K} = \mathcal{H}_1 \oplus \mathcal{H}_2$ to get:

$$\mathcal{K} = \boxed{\begin{array}{c} \mathrm{def}_+(L) = \mathrm{Eig}_{\mathfrak{i}}(L^\star) \xrightarrow{\quad Q \quad} \mathrm{def}_-(L) = \mathrm{Eig}_{-\mathfrak{i}}(L^\star) \\ \mathrm{dom}\, C(L) = (\mathfrak{i} + L)(\mathrm{dom}\, L) \xrightarrow{\quad C \quad} \mathrm{ran}\, C(L) = (\mathfrak{i} - L)(\mathrm{dom}\, L) \end{array}} = \mathcal{K}$$

The structure underlying Theorem 2.20 is a pair of surjective isomorphisms

$$\Psi_\pm : \mathrm{Eig}_{-1}(A^\star B^\star) \to \mathrm{def}_\pm(L^\star) \quad \text{by} \quad \Psi_\pm(u) = \begin{bmatrix} u \\ \pm \mathfrak{i} B^\star u \end{bmatrix}, \quad \text{for } u \in \mathrm{Eig}_{-1}(A^\star B^\star).$$

Consequently, a partial isometry $Q : \mathrm{def}_+(L) \to \mathrm{def}_-(L)$ induces an operator on $\mathrm{Eig}_{-1}(A^\star B^\star)$ (which we denote by $\tilde{Q}$), as the unique operator making the following diagram commute:

$$\begin{array}{ccc} \mathrm{Eig}_{-1}(A^\star B^\star) & \xrightarrow{\tilde{Q}} & \mathrm{Eig}_{-1}(A^\star B^\star) \\ \Psi_+ \downarrow & & \downarrow \Psi_- \\ \mathrm{def}_+(L) & \xrightarrow{Q} & \mathrm{def}_-(L) \end{array}$$

As in the proof of Theorem 2.20(b), we use the notation $L_Q$ for the self-adjoint extension of $L$ corresponding to $Q$. Let us also denote $v := \tilde{Q}u$, for $u \in \mathrm{dom}\,\tilde{Q}$. Then by (2.23), the action of $L_Q$ on a generic element of

$$\mathrm{dom}\, L_Q = \left\{ \begin{bmatrix} x \\ y \end{bmatrix} + \begin{bmatrix} u \\ -\mathfrak{i} B^\star u \end{bmatrix} + Q \begin{bmatrix} u \\ -\mathfrak{i} B^\star u \end{bmatrix} : \begin{bmatrix} x \\ y \end{bmatrix} \in \mathrm{dom}\, L, \begin{bmatrix} u \\ -\mathfrak{i} B^\star u \end{bmatrix} \in \mathrm{def}_+(L) \right\}$$

is given by

$$L_Q\left(\begin{bmatrix} x \\ y \end{bmatrix} + \begin{bmatrix} u \\ -\mathbbm{i}B^\star u \end{bmatrix} + \begin{bmatrix} v \\ \mathbbm{i}B^\star v \end{bmatrix}\right) = \begin{bmatrix} By + \mathbbm{i}u - \mathbbm{i}v \\ Ax + B^\star u + B^\star v \end{bmatrix}$$

where we have used the fact that $L \subseteq L_Q$ to compute the action of $L_Q$ on $x \oplus y$, and (2.21) to compute the action of $L_Q$ on the rest of the expression. Note that since $Q$ is a partial isometry, the corresponding restriction on $\tilde{Q}$ is that

$$\|u\|_1^2 + \|B^\star u\|_2^2 = \left\|\begin{bmatrix} u \\ -\mathbbm{i}B^\star u \end{bmatrix}\right\|_{\mathcal{K}}^2 = \left\|\begin{bmatrix} \tilde{Q}u \\ -\mathbbm{i}B^\star \tilde{Q}u \end{bmatrix}\right\|_{\mathcal{K}}^2 = \|\tilde{Q}u\|_1^2 + \|B^\star \tilde{Q}u\|_2^2.$$

## 3. Gaussian fields and the Malliavin derivative

The term "stochastic calculus" refers loosely to an infinite-dimensional theory of integration/differentiation in the context of stochastic processes; in particular, Itō calculus and its variant Malliavin calculus. The Malliavin calculus is a stochastic version of calculus of variations and provides a robust definition of the derivative of a random variable which allows, e.g., for integration by parts. Below, we study the Malliavin derivative and its adjoint, the Skorokhod integral (an extension of the Itō integral to processes which may not be adapted; cf. [Gaw]). The Skorokhod integral may be interpreted as an infinite-dimensional generalization of the divergence operator. As an example application of symmetric pairs, we show how they can be used to provide a streamlined proof that the Malliavin derivative and Skorokhod integral are closable operators, and in fact are mutually adjoint. The context for this section is Gaussian fields, i.e., Gaussian processes indexed by a Hilbert space $\mathcal{L}$. For example, in the special case where the index Hilbert space is chosen to be $\mathcal{L} = L^2[0, \infty)$, the resulting Gaussian process is Brownian motion, and then the sample points in $\Omega$ may be understood as paths/trajectories. See [Hid, Bel] for background material and further details.

Consider a probability measure space $(\Omega, \Sigma, \mathbb{P})$. We are interested in operators in the Hilbert spaces of random variables

$$\mathcal{H}_1 := L^2(\Omega, \mathbb{P}) \quad \text{and} \quad \mathcal{H}_2 := L^2(\Omega, \mathbb{P}) \otimes \mathcal{L}. \tag{3.1}$$

Note that an element of $\mathcal{H}_1$ is a random variable taking values in $\mathbb{R}$, and $\mathcal{H}_2$ is a random variable taking values in $\mathcal{L}$. We will see in Definition 3.6 that the Malliavin derivative is a densely defined operator on $\mathcal{H}_1$ with range in $\mathcal{H}_2$.

It will be enough for us to work with real Hilbert spaces, in which case we refer to (2.1). We use the notation

$$\mathbb{E}[F] = \int_\Omega F \, d\mathbb{P} \quad \text{and} \quad \mathbb{E}[FG] = \int_\Omega FG \, d\mathbb{P} = \langle F, G \rangle_{\mathcal{H}_1}, \quad \text{for all } F, G \in \mathcal{H}_1.$$

Since $\mathcal{H}_2$ is a tensor product of Hilbert spaces, it inherits a natural inner product from the tensor structure:

$$\langle F_1 \otimes k_1, F_2 \otimes k_2 \rangle_{\mathcal{H}_2} = \langle F_1, F_2 \rangle_{\mathcal{H}_1} \langle k_1, k_2 \rangle_{\mathcal{L}} = \mathbb{E}[F_1 F_2] \langle k_1, k_2 \rangle_{\mathcal{L}}.$$

Suppose that for $i = 1, 2$, we have $\psi_i \in L^2(\Omega, \mathbb{P}) \otimes \mathcal{L}$. In other words, suppose $\psi_i : \Omega \to \mathcal{L}$ are measurable functions on $\Omega$ with

$$\int_\Omega \|\psi_i(\omega)\|^2 \, d\mathbb{P}(\omega) < \infty.$$

Then

$$\langle \psi_1, \psi_2 \rangle_{\mathcal{H}_2} = \int_\Omega \langle \psi_1(\omega), \psi_2(\omega) \rangle_{\mathcal{L}} \, d\mathbb{P}(\omega) = \mathbb{E}[\langle \psi_1(\cdot), \psi_2(\cdot) \rangle_{\mathcal{L}}].$$

**Definition 3.1.** A *Gaussian field* or *Gaussian Hilbert space* is $(\Omega, \Sigma, \mathbb{P}, \Phi)$, where

$$\Phi : \mathcal{L} \to L^2(\Omega, \mathbb{P}),$$

is interpreted as a Wiener process indexed by $\mathcal{L}$ satisfying the following properties:

(i) $\mathbb{E}[\Phi(h)] = 0$, for all $h \in \mathcal{L}$, and

(ii) for any finite subset $\{h_1, \ldots, h_n\} \subseteq \mathcal{L}$, the random variables $\{\Phi(h_i)\}_{i=1}^n$ are jointly Gaussian. That is, density of the joint distribution of $\{\Phi(h_i)\}_{i=1}^n$ is given by
$$\mathfrak{g}_{G_n}(x) := (2\pi)^{-n/2}(\det G_n)^{-1/2}\, e^{-\frac{1}{2}\langle x, G_n^{-1} x\rangle}, \qquad x \in \mathbb{R}^n,$$
where $G_n$ is the Gramian matrix with entries $[G_n]_{i,j} = \langle h_i, h_j\rangle_{\mathcal{L}}$.

**Theorem 3.2.** *Given any real Hilbert space $\mathcal{H}$, there always exist associated Gaussian fields, i.e., an associated probability space and Gaussian process $\Phi$, with properties as specified in Definition 3.1*

A proof of Theorem 3.2 may be found in [Hid] or [PS].

*Remark* 3.3. In the special case of the Ito-Wiener integral
$$\Phi(h) = \int_0^\infty h(t)\, d\Phi_t, \qquad h \in \mathcal{L},$$
we recover from Definition 3.1 the usual Brownian motion, as described in [Hid].

**Definition 3.4.** The *symmetric Fock space* is Hilbert completion of
$$\Gamma_{\text{sym}}(\mathcal{L}) := \bigoplus\nolimits_{n=0}^\infty \text{Sym}(\mathcal{L}^{\otimes n}), \quad \text{where} \quad \mathcal{L}^{\otimes n} := \bigotimes\nolimits_{j=1}^n \mathcal{L}, \tag{3.2}$$
and "Sym" denotes the operator that symmetrizes a tensor. See [BR] for details.

The symmetric Fock space provides a useful model for $\mathcal{H}_1$, as is seen from Lemma 3.5; it is also used in Lemma 3.7 to give insight on the domain of the Malliavin derivative.

**Lemma 3.5.** *The symmetric Fock space is isometrically isomorphic to $\mathcal{H}_1 = L^2(\Omega, \Sigma, \mathbb{P})$.*

*Proof.* For $k \in \mathcal{L}$, define $e^k := \sum_{n=0}^\infty \frac{k^{\otimes n}}{\sqrt{n!}}$. Then for $k_i \in \mathcal{L}$, $i = 1, 2$, the vectors $e^{k_i}$ are the unique elements of $\Gamma_{\text{sym}}(\mathcal{L})$ such that
$$\langle e^{k_1}, e^{k_2}\rangle_{\Gamma_{\text{sym}}(\mathcal{L})} = \sum_{n=0}^\infty \frac{\langle k_1, k_2\rangle_{\mathcal{L}}^n}{n!} = e^{\langle k_1, k_2\rangle_{\mathcal{L}}}. \tag{3.3}$$
Moreover, the mapping
$$W_0 : \{e^k\}_{k \in \mathcal{L}} \to L^2(\Omega, \mathbb{P}) \qquad \text{by} \qquad W_0(e^k) = e^{\Phi(k) - \frac{1}{2}\|k\|_2^2} \tag{3.4}$$
extends by linearity and closure to a unitary isomorphism $W : \Gamma_{\text{sym}}(\mathcal{L}) \to L^2(\Omega, \mathbb{P})$. □

**Definition 3.6.** For $n \in \mathbb{N}$, let $p \in \mathbb{R}[x_1, \ldots, x_n]$ be a polynomial in $n$ real variables, so $p : \mathbb{R}^n \to \mathbb{R}$. Define the *Malliavin derivative* by
$$\text{dom}\, T := \text{span}_{p \in \mathbb{R}[x_1, \ldots, x_n]}\{p(\Phi(h_1), \ldots, \Phi(h_n)) : (h_i)_{i=1}^n \subseteq \mathcal{L}, \text{ and } n \in \mathbb{N}\}, \tag{3.5}$$
$$T(F) := \sum_{j=1}^n \left(\frac{\partial p}{\partial x_j}\right)(\Phi(h_1), \ldots, \Phi(h_n)) \otimes h_j, \quad \text{for } F \in \text{dom}\, T. \tag{3.6}$$
Note that and $T(F) \in \mathcal{H}_2$, as in (3.1).

**Lemma 3.7.** dom $T$ *is dense in $\mathcal{H}_1$.*

*Proof.* Since span$\{e^k : k \in \mathcal{L}\}$ is dense in $\Gamma_{\text{sym}}(\mathcal{L})$, it follows that dom $T$ is dense in $L^2(\Omega, \mathbb{P})$ by Lemma 3.5. □

**Definition 3.8.** For dom $T$ as in (3.5), define the stochastic integral operator
$$\text{dom}\, S := \text{dom}\, T \otimes \mathcal{L} \tag{3.7}$$
$$S(F \otimes k) := F \cdot \Phi(k) - \langle T(F), k\rangle. \tag{3.8}$$
Note that (3.7) is an algebraic tensor product, not a Hilbert tensor product. It is immediate that $S$ is densely defined because the span of the elementary tensors is dense in the tensor product.

*Remark* 3.9. We will see in Theorem 2.5 below that $(S, T)$ is a symmetric pair in the sense of Definition 2.1.

$$\mathcal{H}_1 = L^2(\Omega, \mathbb{P}) \xrightarrow{T} \mathcal{H}_2 = L^2(\Omega, \mathbb{P}) \otimes \mathcal{L}$$
$$\xleftarrow{S}$$

**Definition 3.10.** We denote the standard Gaussian density by

$$\mathfrak{g}_n(x) = (2\pi)^{-n/2} \, e^{-\frac{1}{2}\|x\|^2}, \qquad x \in \mathbb{R}^n, \tag{3.9}$$

and recall that the corresponding expectation is given by

$$\mathbb{E}\left[p\left(\Phi(h_1), \ldots, \Phi(h_n)\right)\right] := \int_{\mathbb{R}^n} p(x_1, \ldots, x_n) \mathfrak{g}_n(x_1, \ldots, x_n) \, dx_1 \ldots dx_n \tag{3.10}$$

where the vectors $\{h_i\}_{i=1}^n \subseteq \mathcal{L}$ are chosen to be orthogonal, i.e., with $\langle h_i, h_j \rangle_{\mathcal{H}} = \delta_{ij}$. Without this restriction, the joint distribution of the vector-valued random variable $(\Phi(h_1), \ldots, \Phi(h_n)) \in \mathbb{R}^n$ is the $n$-dimensional Gaussian whose covariance matrix is the Gramian with entries $\langle h_i, h_j \rangle_{\mathcal{H}}$.

**Lemma 3.11.** *For all $k \in \mathcal{L}$ and $F \in \operatorname{dom} T$, we have*

$$\langle T(F), \mathbf{1} \otimes k \rangle_{\mathcal{L}} \rangle_{\mathcal{H}_2} = \langle F, \Phi(k) \rangle_{\mathcal{H}_1}. \tag{3.11}$$

*Proof.* Note that the tensor structure of $\mathcal{H}_2 = L^2(\Omega, \mathbb{P}) \otimes \mathcal{L}$ means

$$\langle T(F), \mathbf{1} \otimes k \rangle_{\mathcal{L}} \rangle_{\mathcal{H}_2} = \mathbb{E}\left[\langle T(F), k \rangle_{\mathcal{L}}\right], \tag{3.12}$$

where the expectation on the right is with respect to the Gaussian density as in (3.10). Consequently, we may prove (3.11) in the form

$$\mathbb{E}[\langle T(F), k \rangle_{\mathcal{L}}] = \mathbb{E}[F \cdot \Phi(k)], \qquad \text{for all } k \in \mathcal{L}, F \in \operatorname{dom} T. \tag{3.13}$$

Choosing an $F$ includes a choice of $(\ell_i)_{i=1}^n$. The calculations which follow will be simplify greatly when $k$ is orthogonal to $\ell_i$. We can effect this by applying the Gram-Schmidt procedure to this set, which will only result in a different polynomial in the representation of $T$ (and we are free to range over all polynomials). Thus we may pass from $(\ell_i)_{i=1}^n$ to $(h_i)_{i=1}^n$, where $\langle h_i, h_j \rangle_{\mathcal{L}} = \delta_{ij}$. In other words, applying the spectral theorem to the Gramian of $(\ell_i)_{i=1}^n$ allows us to diagonalize this matrix so that we can work instead with the standard Gaussian covariance matrix $\mathbb{I}$. By choosing $h_1 = k$ in the Gram-Schmidt procedure, we have $\langle k, h_i \rangle_{\mathcal{L}} = 1$ for $i = 1$ and 0 otherwise. Consequently, applying orthogonality and integration by parts to (3.6) yields

$$\begin{aligned}
\mathbb{E}\left[\langle T(F), k \rangle_{\mathcal{L}}\right] &= \mathbb{E}\left[\frac{\partial p}{\partial x_1}\left(\Phi(h_1), \ldots, \Phi(h_n)\right)\right] & \langle k, h_i \rangle_{\mathcal{L}} = \delta_{1i} \\
&= \int_{\mathbb{R}^n} \frac{\partial p}{\partial x_1}(x_1, \ldots, x_n) \mathfrak{g}_n(x_1, \ldots, x_n) \, dx_1 \ldots dx_n & \text{by (3.10)} \\
&= -\int_{\mathbb{R}^n} p(x_1, \ldots, x_n) \frac{\partial \mathfrak{g}_n}{\partial x_1}(x_1, \ldots, x_n) \, dx_1 \ldots dx_n & \text{IBP} \\
&= \int_{\mathbb{R}^n} x_1 p(x_1, \ldots, x_n) \mathfrak{g}_n(x_1, \ldots, x_n) \, dx_1 \ldots dx_n & \text{by (3.9)} \\
&= \mathbb{E}\left[p\left(\Phi(h_1), \ldots, \Phi(h_n)\right) \Phi(h_1)\right] & \text{by (3.10)} \\
&= \mathbb{E}[F \cdot \Phi(k)] & \Phi(k) = \Phi(h_1),
\end{aligned}$$

which is (3.13). $\square$

While using Lemma 3.11 to derive our main result Theorem 2.5, we will need the following simple result about the Leibnitz rule in this context.

**Lemma 3.12.** *The Malliavin derivative $T$ is a module derivation. In other words,*

$$T(HK) = T(H)K + HT(K), \qquad \text{for any } H, K \in \operatorname{dom} T. \tag{3.14}$$

*Proof.* Since the vector space of all polynomials of degree at most $n$ is invariant under changes of coordinates, it suffices to note that we can find $n$ and $h_1, \ldots, h_n \in \mathcal{L}$ which can be used to represent both $H$ and $K$ simultaneously:

$$H = p(\Phi(h_1), \ldots, \Phi(h_n)) \quad \text{and} \quad K = q(\Phi(h_1), \ldots, \Phi(h_n)).$$

Now (3.14) follows from the fact that $\frac{\partial (pq)}{\partial x_i} = \frac{\partial p}{\partial x_i} q + p \frac{\partial q}{\partial x_i}$. $\square$

**Theorem 3.13.** *The Malliavin derivative $T$ of Definition 3.6 and the stochastic integral operator $S$ of Definition 3.8 form a symmetric pair:*

$$\langle T(F), G \otimes k \rangle_{\mathcal{H}_2} = \langle F, S(G \otimes k) \rangle_{\mathcal{H}_1}, \qquad \text{for any } F, G \in \operatorname{dom} T, k \in \mathcal{L}. \tag{3.15}$$

*Consequently, $T$ and $S$ are closable linear operators. Furthermore, $(T, S)$ is a maximal symmetric pair, whence $T^\star = \overline{S}$ and $S^\star = \overline{T}$.*

*Proof.* We compute (3.15) first:

$$\begin{aligned}
\langle T(F), G \otimes k \rangle_{\mathcal{H}_2} &= \mathbb{E}\left[ \langle T(F), k \rangle_{\mathcal{L}} G \right] & \text{by (3.12)} \\
&= \mathbb{E}\left[ \langle T(FG), k \rangle_{\mathcal{L}} \right] - \mathbb{E}\left[ F \langle T(G), k \rangle_{\mathcal{L}} \right] & \text{by (3.14)} \\
&= \mathbb{E}\left[ FG \cdot \Phi(k) \right] - \mathbb{E}\left[ F \langle T(G), k \rangle_{\mathcal{L}} \right] & \text{by (3.13)} \\
&= \mathbb{E}\left[ F \left( G \cdot \Phi(k) - \langle T(G), k \rangle_{\mathcal{L}} \right) \right] & \text{by linearity} \\
&= \mathbb{E}\left[ F \cdot S(G \otimes k) \right] & \text{by (3.8)} \\
&= \langle F, S(G \otimes k) \rangle_{\mathcal{H}_1}.
\end{aligned}$$

The closability of $T$ and $S$ now follows from Theorem 2.5.

To see maximality, we show that no element of $\Gamma(S^\star)$ is orthogonal to $\Gamma(\overline{T})$. More precisely, since $\overline{T}(\mathrm{e}^{\Phi(k)}) = \mathrm{e}^{\Phi(k)} \otimes k$ by Corollary 3.17, just below, we show that if $F \in \operatorname{dom} S^\star$ and

$$\left\langle \begin{bmatrix} \mathrm{e}^{\Phi(k)} \\ \mathrm{e}^{\Phi(k)} \otimes k \end{bmatrix}, \begin{bmatrix} F \\ S^\star(F) \end{bmatrix} \right\rangle = 0, \qquad \text{for all } k \in \mathcal{L}, \tag{3.16}$$

then $F = 0$. Observe that (3.16) is equivalent to

$$\mathbb{E}\left[ \mathrm{e}^{\Phi(k)} (F + \langle S^\star(F), k \rangle_{\mathcal{L}}) \right] = 0, \qquad \text{for all } k \in \mathcal{L}, \tag{3.17}$$

and since $\operatorname{span}\{\mathrm{e}^{\Phi(k)}\}_{k \in \mathcal{L}}$ is dense in $\mathcal{H}_1$, (3.17) implies $F = 0$. $\square$

**Corollary 3.14.** *For $T_k(F) := \langle T(F), k \rangle_{\mathcal{L}}$, we have that $T_k$ is a derivation and that*

$$T_k + T_k^\star = M_{\Phi(k)},$$

*where $M_f$ denotes the multiplication operator (by the function $f$).*

*Proof.* By (3.8), we have $S(F \otimes k) = M_{\Phi(k)} F - T_k(F)$, so that Theorem 2.5 gives

$$(T_k + T_k^\star)(F) = T_k(F) + M_{\Phi(k)} F - T_k(F) = M_{\Phi(k)} F, \qquad \text{for all } F \in \operatorname{dom} T.$$

It is immediate from (3.6) that $T_k$ is a derivation. $\square$

*Remark* 3.15. The operator $S$ of Definition 3.8 is essentially integration by parts, and this is highlighted by the proof of Theorem 2.5, which is essentially the same as the calculus proof of integration-by-parts based on the product rule. Example 3.16 extends this parallel.

**Example 3.16** (Vector bundles)**.** The construction above is an infinite-dimensional generalization of vector bundles. For a vector field $X$ on a manifold $M$, and $f \in C^\infty(M)$, let

$$fX(\varphi)(m) = f(m) X(\varphi)(m), \qquad \forall m \in M, \varphi \in C^\infty(M).$$

Then $(fX)^\star = -fX + M_{X(f)}$.

**Corollary 3.17.** *Let $(\Omega, \Sigma, \mathbb{P}, \Phi)$ be a Gaussian field and let $\overline{T} : \mathcal{H}_1 \to \mathcal{H}_2$ be the closure of $T$ as in Definition 3.6. For $k \in \mathcal{L}$, the exponential random variables $\mathrm{e}^{\Phi(k)}$ lie in $\operatorname{dom}(T^\star \overline{T})$, and*

$$T^\star \overline{T}(\mathrm{e}^{\Phi(k)}) = \left( \Phi(k) - \|k\|_{\mathcal{L}}^2 \right) \mathrm{e}^{\Phi(k)}.$$

See Remark 2.8 for a discussion of the significance of the following result on $\ker(\overline{T})$.

**Corollary 3.18.** *Let $\mathcal{T}$, $\mathcal{L}$, and $(\Omega, \Sigma, \mathbb{P}, \Phi)$ be as above. Then $\ker(\overline{T})$ is the 1-dimensional subspace of $\mathcal{H}_1 = L^2(\Omega, \mathbb{P})$ spanned by the constant function $\mathbf{1}$ on $\Omega$.*

*Proof.* Let $p_n^H$ denote the Hermite polynomial of degree $n$, so that

$$e^{tx - \frac{t^2}{2}} = \sum_{n=0}^{\infty} p_n^H(x) \frac{t^n}{n!}$$

By (3.4), the isomorphism $W : \Gamma_{\text{sym}}(\mathcal{L}) \to L^2(\Omega, \mathbb{P})$ of Lemma 3.7 maps $H_n$ (as defined in (3.2)) into $\text{span}\{p_n^H(\Phi(k))\}_{k \in \mathcal{L}}$. For $F \in \text{dom}\,\overline{T}$ with $\overline{T}(F) = 0$ one has $T^\star \overline{T}(F) = 0$, so that

$$0 = T^\star T(F) = \sum_{n=0}^{\infty} n F_n,$$

where $F_n \in W(H_n)$. Therefore, $n F_n = 0$ for $n = 1, 2, \ldots$, and hence $F_n = 0$ for $n = 1, 2, \ldots$. It follows that $F_n \neq 0$ can only hold for $n = 0$. Since $W(H_0) = \text{span}\,\mathbf{1}$, the conclusion follows. □

## 4. Tomita-Takesaki theory

Part of von Neumann's classification of rings of operators (now called von Neumann algebras) included the classification of factors; of these the type III factors were considered intractable up to the 1970s. The primary reason for this, is the non-availability of even semifinite traces; see the discussion in §1.2. The Tomita-Takesaki theory was a response to this obstacle, and it has proved successful, leading to a structure theory for these previously intractable cases. The part of the theory of relevance to us here is the modular conjugation (and automorphisms), and its connection to a certain closable operator, as well as its associated polar decomposition. See [Kad, Tak1, Tak2, KR1, BR, HV] for further background.

Let $\mathcal{H}$ be a Hilbert space and consider the bounded linear operators on $\mathcal{H}$, denoted $B(\mathcal{H})$. We use the usual bracket notation to denote the commutator of two operators:

$$[x, y] := xy - yx, \qquad x, y \in B(\mathcal{H}).$$

**Definition 4.1.** For a von Neumann algebra $\mathfrak{M} \subseteq B(\mathcal{H})$, the *commutator* of $\mathfrak{M}$ is

$$\mathfrak{M}' := \{x \in B(\mathcal{H}) : [x, m] = 0, \forall m \in \mathfrak{M}\}. \tag{4.1}$$

**Definition 4.2.** For $\xi \in \mathcal{H}$ with $\|\xi\| = 1$, we say that $\xi$ is $\mathfrak{M}$-*cyclic* iff

$$\mathfrak{M}\xi := \{m\xi : m \in \mathfrak{M}\} \tag{4.2}$$

is dense in $\mathcal{H}$, and we say that $\xi$ is $\mathfrak{M}$-*separating* iff

$$\mathfrak{M}'\xi := \{m'\xi : m' \in \mathfrak{M}'\} \tag{4.3}$$

is dense in $\mathcal{H}$.

*Remark* 4.3. Throughout the rest of this section, we assume that a cyclic and separating vector $\xi \in \mathcal{H}$ has been fixed. It can be shown that $\mathfrak{M}'\xi$ is dense in $\mathcal{H}$ iff

$$m\xi = 0,\ m \in \mathfrak{M} \quad \implies \quad m = 0, \tag{4.4}$$

and this is the motivation for the terminology *separating* in Definition 4.2.

Given $\mathfrak{M}$ and a fixed faithful normal state, one may construct a pair $(\mathcal{H}, \xi)$ via the GNS construction such that $\xi$ will be a cyclic and separating vector in $\mathcal{H}$ for $\mathfrak{M}$. See [KR2, BR, vD] for background material and further detail. However, in general a von Neumann algebra $\mathfrak{M} \subseteq B(\mathcal{H})$ with $\dim \mathcal{H} = \aleph_0$ cannot be realized with cyclic and separating vectors. For example, with $\mathfrak{M} = B(\mathcal{H})$ and $\mathfrak{M}' = \mathbb{C}\mathbb{I}_\mathcal{H}$, there is no $\xi \in \mathcal{H}$ for which $\mathfrak{M}'\xi$ is dense in $\mathcal{H}$. We mention this just to emphasize the importance of the assumption of a cyclic and separating vector $\xi \in \mathcal{H}$.

**Definition 4.4.** Define the (unbounded) conjugate-linear operators $S_\xi$ and $F_\xi$ by

$$\text{dom}\, S_\xi := \mathfrak{M}\xi, \qquad S_\xi(m\xi) := m^\star \xi, \quad \text{for all } m\xi \in \mathfrak{M}\xi. \tag{4.5}$$

and

$$\text{dom}\, F_\xi := \mathfrak{M}'\xi, \qquad F_\xi(m'\xi) := (m')^\star \xi, \quad \text{for all } m'\xi \in \mathfrak{M}'\xi. \tag{4.6}$$

Note that by our choice of $\xi$, both $S_\xi$ and $F_\xi$ are densely defined; see Remark 4.3 and Definition 4.2.

We will show that $(S_\xi, F_\xi)$ is a symmetric pair as follows:

$$\mathfrak{M}\xi \subseteq \mathcal{H} \underset{F_\xi}{\overset{S_\xi}{\rightleftarrows}} \mathfrak{M}'\xi \subseteq \mathcal{H}$$

Note that since $S_\xi$ and $F_\xi$ are conjugate-linear operators, we refer to (2.3).

*Remark* 4.5. It is immediate from (4.2), (4.2), and (4.4) that $\ker S_\xi = 0$ and $\overline{\text{ran}\, S_\xi} = \mathcal{H}$.

**Theorem 4.6.** $(S_\xi, F_\xi)$ *is a symmetric pair, and thus* $S_\xi$ *and* $F_\xi$ *are closable operators and*

$$S_\xi \subseteq F_\xi^\star \qquad \text{and} \qquad F_\xi \subseteq S_\xi^\star. \tag{4.7}$$

*Proof.* Take $u = m\xi \in \mathfrak{M}\xi$ and $v = m'\xi \in \mathfrak{M}'\xi$. Then (4.5) gives

$$\langle S_\xi u, v\rangle_\mathcal{H} = \langle S_\xi m\xi, m'\xi\rangle_\mathcal{H} = \langle m^\star\xi, m'\xi\rangle_\mathcal{H} = \langle \xi, mm'\xi\rangle_\mathcal{H}, = \langle \xi, m'm\xi\rangle_\mathcal{H},$$

since $m' \in \mathfrak{M}'$ implies $mm' = m'm$. Then

$$\langle S_\xi u, v\rangle_\mathcal{H} = \langle \xi, m'm\xi\rangle_\mathcal{H} = \overline{\langle \xi, (m'm)^\star\xi\rangle}_\mathcal{H} = \overline{\langle \xi, m^\star(m')^\star\xi\rangle}_\mathcal{H} = \overline{\langle m\xi, F_\xi(m'\xi)\rangle}_\mathcal{H} = \overline{\langle u, F_\xi v\rangle}_\mathcal{H},$$

so $(S_\xi, F_\xi)$ is a symmetric pair. Theorem 2.5 ensures $S_\xi$ and $F_\xi$ are closable. □

*Remark* 4.7. The significance of Theorem 4.6 is that it allows for the polar decomposition of $\overline{S_\xi}$ (and $\overline{F_\xi}$), as in Definition 4.8, which is essential for Tomita-Takesaki theory. See [DS] for details on the polar decomposition in the unbounded case.

**Definition 4.8.** In view of Theorem 4.6, the *modular operator*

$$\Delta : \mathcal{H} \to \mathcal{H} \qquad \text{by} \qquad \Delta := (S_\xi^\star \overline{S}_\xi)^{1/2} \tag{4.8}$$

is well-defined; it is well-known to be self-adjoint and positive. Also, the *modular conjugation* is the conjugate linear partial isometry

$$J : \mathcal{H} \to \mathcal{H} \qquad \text{defined by} \qquad \overline{S}_\xi = J\Delta. \tag{4.9}$$

*Remark* 4.9. As mentioned in §1.2, the foundational results of Tomita-Takesaki are

$$J\mathfrak{M}J = \mathfrak{M}' \qquad \text{and} \qquad \Delta^{\text{i}t}\mathfrak{M}\Delta^{-\text{i}t} = \mathfrak{M}, \quad \text{for all } t \in \mathbb{R}. \tag{4.10}$$

These results are used to develop a robust structure theory for von Neumann algebras which are type III factors. A *factor* is a von Neumann algebra whose center consists only of multiples of the identity operator, and a factor is said to be *type III* iff $\mathfrak{M}$ does not have any (even semifinite) trace; cf. §1.2. However, as a precursor to (4.10), it must be shown that $S_\xi$ is closable, so that (4.9) makes sense. See [BR, HV] for details.

**Theorem 4.10.** *With $J$ as defined in* (4.9), $J\mathfrak{M}J = \mathfrak{M}'$.

*Sketch of proof.* On $\text{dom}\, S_\xi$, we have the identity

$$(S_\xi x S_\xi) y = y(S_\xi x S_\xi), \qquad \text{for all } x, y \in \mathfrak{M}, \tag{4.11}$$

which can be verified directly from the definition of $S_\xi$. Note that the operators $x, y \in \mathfrak{M}$ appearing in (4.11) are understood to be acting on their natural dense domains in $\mathcal{H}$ so that each side of the equality makes sense; see the last part of Remark 4.3. Equation (4.11) indicates that for any $x \in \mathfrak{M}$, one has $S_\xi x S_\xi \in \mathfrak{M}'$. Now (4.8)–(4.9) and von Neumann's double commutant theorem $\mathfrak{M}'' = \mathfrak{M}$ imply that

$$JxJ \in \mathfrak{M}', \qquad \text{for all } x \in \mathfrak{M}. \tag{4.12}$$

It follows from $S_\xi^2 = \mathbb{I}$, that $J^\star J = \mathbb{I}$, and thus (4.12) implies $J\mathfrak{M}J = \mathfrak{M}'$. □

**Theorem 4.11.** *The pair $(S_\xi, F_\xi)$ is a maximal symmetric pair, so that $\overline{S_\xi} = F_\xi^\star$ and $\overline{F_\xi} = S_\xi^\star$.*

*Proof.* From Theorem 4.6, we have the graph containment $\Gamma(S_\xi) \subseteq \Gamma(F_\xi^\star)$. To prove that this is actually an equality, we show that

$$\begin{bmatrix} \zeta \\ F_\xi^\star \zeta \end{bmatrix} \in \Gamma(F_\xi^\star) \text{ and } \begin{bmatrix} \zeta \\ F_\xi^\star \zeta \end{bmatrix} \perp \Gamma(S_\xi) \qquad \implies \qquad \zeta = 0. \tag{4.13}$$

From the hypothesis in (4.13), we have the identity

$$\langle m\xi, \zeta \rangle_{\mathcal{H}} + \langle m^\star \xi, F_\xi^\star \zeta \rangle_{\mathcal{H}} = 0, \qquad \text{for all } m \in \mathfrak{M}. \tag{4.14}$$

If we take $m = m^\star \in \mathfrak{M}$ in (4.14), it implies

$$F_\xi^\star \zeta = -\zeta \tag{4.15}$$

because $\xi$ is cyclic. By way of contradicting (4.13), suppose that $\zeta \neq 0$. Then for $\varepsilon_k$ satisfying $0 < \varepsilon_k < \frac{1}{2}\|\zeta\|^2$, we can choose $m'_k \in \mathfrak{M}'$ satisfying

$$\text{(i) } m'_k = (m'_k)^\star \qquad \text{and} \qquad \text{(ii) } \left| \langle m'_k \xi, \zeta \rangle - \|\zeta\|^2 \right| < \varepsilon_k \tag{4.16}$$

because $\xi$ is $\mathfrak{M}$-separating, and hence $\mathfrak{M}'$-cyclic. But then

$$\begin{aligned}
\langle m'_k \xi, \zeta \rangle_{\mathcal{H}} &= -\langle m'_k \xi, F_\xi^\star \zeta \rangle_{\mathcal{H}} && \text{by (4.15)} \\
&= -\langle (m'_k)^\star \xi, F_\xi^\star \zeta \rangle_{\mathcal{H}} && \text{by (4.16i)} \\
&= -\langle (m'_k)^\star \xi, \zeta \rangle_{\mathcal{H}} && \text{by (4.6)} \\
&= -\langle m'_k \xi, \zeta \rangle_{\mathcal{H}} && \text{by (4.16i) again.}
\end{aligned}$$

Now (4.16ii) implies $\|\zeta\|^2 = -\|\zeta\|^2$, whence $\zeta = 0$. □

## 5. Laplace operators on infinite networks

In this section, we discuss two Hilbert spaces of functions defined on certain graphs (called resistance networks; see Definition 5.2). One is the familiar $\ell$ space of square-summable functions, which we denote $\ell^2(G)$, where $G$ is a resistance network. The other Hilbert space is the slightly less familiar space $\mathcal{H}_\mathcal{E}$ of functions of finite energy (or Dirichlet functions). These Hilbert spaces have different inner products, and it is interesting to compare the Laplace operator and its spectral properties in these two different contexts. For further background on resistance networks, the energy form $\mathcal{E}$, the Laplace operator $\Delta$, and other related topics in discrete resistance analysis, we refer the reader to [JP7, JP3, JP5, JP4, JP6, JP8, JP9, JP2, JP1, LP, KL1, KL2].

*Remark* 5.1. The Laplace operator $\Delta$ studied in this section is not related to the modular operator of the previous section, which is unfortunately traditionally also denoted by the same symbol.

We revisit the example presented in [JP10] and show that the Laplacian can be considered as an operator from $\mathcal{H}_\mathcal{E}$ to $\ell^2(G)$. If we define the "inclusion" operator $K : \ell^2(G) \to \mathcal{H}_\mathcal{E}$ by $K\varphi = \varphi$ on the appropriate dense domain in $\ell^2(G)$, then $(K, \Delta)$ forms a symmetric pair. In contrast to the examples in §3 and §4, the construction here illustrates a symmetric pair which is *not* maximal; see Remark 5.19.

**Definition 5.2.** A *(resistance) network* $(G, c)$ is a connected weighted undirected graph with vertex set $G$ and adjacency relation defined by a symmetric *conductance function* $c : G \times G \to [0, \infty)$. More precisely, there is an edge connecting $x$ and $y$ iff $c_{xy} > 0$, in which case we write $x \sim y$. The nonnegative number $c_{xy} = c_{yx}$ is the weight associated to this edge and it is interpreted as the conductance, or reciprocal resistance of the edge.

We make the standing assumption that $(G, c)$ is *locally finite*. This means that every vertex has *finite degree*, i.e., for any fixed $x \in G$ there are only finitely many $y \in G$ for which $c_{xy} > 0$. We denote the net conductance at a vertex by

$$c(x) := \sum_{y \sim x} c_{xy}. \tag{5.1}$$

Motivated by current flow in electrical networks, we also assume $c_{xx} = 0$ for every vertex $x \in G$.

In this paper, *connected* means simply that for any $x, y \in G$, there is a finite sequence $\{x_i\}_{i=0}^n$ with $x = x_0$, $y = x_n$, and $c_{x_{i-1}x_i} > 0$, $i = 1, \ldots, n$. For any network, one can fix a reference vertex, which we shall denote by $o$ (for "origin"). It will always be apparent that our calculations depend in no way on the choice of $o$.

**Definition 5.3.** The *energy form* is the (closed, bilinear) Dirichlet form

$$\mathcal{E}(u, v) := \frac{1}{2} \sum_{x,y \in G} c_{xy}(u(x) - u(y))(v(x) - v(y)), \tag{5.2}$$

which is defined whenever the functions $u$ and $v$ lie in the domain

$$\operatorname{dom} \mathcal{E} = \{u : G \to \mathbb{R} : \mathcal{E}(u, u) < \infty\}. \tag{5.3}$$

Hereafter, we write the energy of $u$ as $\mathcal{E}(u) := \mathcal{E}(u, u)$. Note that $\mathcal{E}(u)$ is a sum of nonnegative terms and hence converges iff it converges absolutely.

The energy form $\mathcal{E}$ is sesquilinear and conjugate symmetric on $\operatorname{dom} \mathcal{E}$ and would be an inner product if it were positive definite. Let $\mathbf{1}$ denote the constant function with value 1 and observe that $\ker \mathcal{E} = \mathbb{R}\mathbf{1}$. One can show that $\operatorname{dom} \mathcal{E}/\mathbb{R}\mathbf{1}$ is complete and that $\mathcal{E}$ is closed; see [JP7, JP1], [Kat], or [FŌT].

**Definition 5.4.** The *energy (Hilbert) space* is $\mathcal{H}_\mathcal{E} := \operatorname{dom} \mathcal{E}/\mathbb{R}\mathbf{1}$. The inner product and corresponding norm are denoted by

$$\langle u, v \rangle_\mathcal{E} := \mathcal{E}(u, v) \quad \text{and} \quad \|u\|_\mathcal{E} := \mathcal{E}(u, u)^{1/2}. \tag{5.4}$$

It is shown in [JP7, Lem. 2.5] that the evaluation functionals $L_x u = u(x) - u(o)$ are continuous, and hence correspond to elements of $\mathcal{H}_\mathcal{E}$ by Riesz duality (see also [JP7, Cor. 2.6]). When considering $\mathbb{C}$-valued functions, (5.4) is modified as follows: $\langle u, v \rangle_\mathcal{E} := \mathcal{E}(\overline{u}, v)$.

**Definition 5.5.** Let $v_x$ be defined to be the unique element of $\mathcal{H}_\mathcal{E}$ for which

$$\langle v_x, u \rangle_\mathcal{E} = u(x) - u(o), \qquad \text{for every } u \in \mathcal{H}_\mathcal{E}. \tag{5.5}$$

Note that $v_o$ corresponds to a constant function, since $\langle v_o, u \rangle_\mathcal{E} = 0$ for every $u \in \mathcal{H}_\mathcal{E}$. Therefore, $v_o$ may be safely omitted in some calculations.

*Remark* 5.6. As (5.5) means that the collection $\{v_x\}_{x \in G}$ forms a reproducing kernel for $\mathcal{H}_\mathcal{E}$, we call $\{v_x\}_{x \in G}$ the *energy kernel*. The energy kernel has dense span in $\mathcal{H}_\mathcal{E}$; cf. [Aro]. To see this, note that a RKHS is a Hilbert space $H$ of functions on some set $X$, such that point evaluation by points in $X$ is continuous in the norm of $H$. Consequently, every $x \in X$ defines a vector $k_x \in H$ by Riesz's Theorem, and it is immediate from this that $\operatorname{span}\{k_x\}_{x \in X}$ is dense in $H$.

**Definition 5.7.** Let $\delta_x \in \ell^2(G)$ denote the Dirac mass at $x$, i.e., the characteristic function of the singleton $\{x\}$ and let $\delta_x \in \mathcal{H}_\mathcal{E}$ denote the element of $\mathcal{H}_\mathcal{E}$ which has $\delta_x \in \ell^2(G)$ as a representative. The context will make it clear which meaning is intended.

*Remark* 5.8. Observe that $\mathcal{E}(\delta_x) = c(x) < \infty$ is immediate from (5.2), and hence one always has $\delta_x \in \mathcal{H}_\mathcal{E}$ (recall that $c(x)$ is the total conductance at $x$; see (5.1)).

**Definition 5.9.** For $v \in \mathcal{H}_\mathcal{E}$, one says that $v$ has *finite support* iff there is a finite set $F \subseteq G$ such that $v(x) = k \in \mathbb{C}$ for all $x \notin F$. Equivalently, the set of functions of finite support is

$$\operatorname{span}\{\delta_x\} = \{u \in \operatorname{dom} \mathcal{E} : u \text{ is constant outside some finite set}\}. \tag{5.6}$$

Define $\mathcal{F}in$ to be the $\mathcal{E}$-closure of $\operatorname{span}\{\delta_x\}$.

**Definition 5.10.** The set of harmonic functions of finite energy is denoted

$$\mathcal{H}arm := \{v \in \mathcal{H}_\mathcal{E} : \Delta v(x) = 0, \text{ for all } x \in G\}. \tag{5.7}$$

**Theorem 5.11** (Royden Decomposition). $\mathcal{H}_\mathcal{E} = \mathcal{F}in \oplus \mathcal{H}arm$.

Theorem 5.11 is well known and first appeared in [Yam, Thm. 4.1], where it was called the "Royden Decomposition" by analogy with Royden's result for Riemann surfaces. However, the result is incorrect as stated there and the present corrected form may be found in [Soa, §VI] or [LP, §9.3]. This result also follows immediately from Lemma 5.15; see [JP7, Thm. 2.15].

**Definition 5.12.** The *Laplacian* on $G$ is the linear difference operator which acts on a function $u : G \to \mathbb{R}$ by
$$(\Delta u)(x) := \sum_{y \sim x} c_{xy}(u(x) - u(y)). \tag{5.8}$$

We will consider $\Delta$ as an operator in two distinct senses:

$$\begin{align}
(i) \quad & \Delta : \mathcal{H}_\mathcal{E} \to \ell^2(G) && \text{with domain} && \operatorname{dom} \Delta := \operatorname{span}\{v_x\}_{x \in G} & (5.9)\\
(ii) \quad & \Delta_\mathcal{E} : \mathcal{H}_\mathcal{E} \to \mathcal{H}_\mathcal{E} && \text{with domain} && \operatorname{dom} \Delta_\mathcal{E} := \operatorname{span}\{v_x\}_{x \in G}. & (5.10)
\end{align}$$

The first case will be used in the discussion of symmetric pairs; see. The second case will be distinguished by the subscript, and the usage coincides with previous papers; cf. [JP7, JP3, JP5, JP4, JP6, JP8, JP9, JP2, JP1].

In both cases, the operator has the same domain (and is thus densely defined by Remark 5.6), and is determined on its domain by formula (5.8). Note that the sum in (5.8) is finite by the local finiteness assumption in Definition 5.2, and so the Laplacian is well-defined.

**Definition 5.13.** A *dipole* is any $v \in \mathcal{H}_\mathcal{E}$ satisfying the pointwise identity $\Delta v = \delta_x - \delta_y$ for some $x, y \in G$.

*Remark* 5.14. It is easy to see from the definitions (or [JP7, Lemma 2.13]) that energy kernel elements are dipoles. Dipoles may or may not be unique; this is determined by the asymptotic growth rate of $G$ and the properties of the conductance function $c$; see [JP7, JP3, JP5, JP4, JP1]. In fact, one could equivalently define $\Delta$ on $\operatorname{dom} \Delta = \operatorname{span}\{v_x\}_{x \in G}$ (in either of the senses in Definition 5.12) by
$$\Delta v_x = \delta_x - \delta_o, \tag{5.11}$$
in place of (5.8). It follows that one can therefore always find a dipole for any given pair of vertices $x, y \in G$, namely, $v_x - v_y$.

**Lemma 5.15** ([JP7, Lem. 2.11]). *For every* $x \in G$,
$$\langle \delta_x, u \rangle_\mathcal{E} = \Delta u(x), \quad \text{for all } u \in \operatorname{dom} \Delta. \tag{5.12}$$
*Similarly,* $\langle \delta_x, u \rangle_\mathcal{E} = \Delta_\mathcal{E} u(x)$, *for all* $u \in \operatorname{dom} \Delta_\mathcal{E}$.

*Proof.* One can compute $\langle \delta_x, u \rangle_\mathcal{E} = \mathcal{E}(\delta_x, u)$ directly from formula (5.2). □

**Definition 5.16.** For functions $u, v : G \to \mathbb{R}$, we also have the usual the inner product of $\ell^2(G)$ defined by
$$\langle u, v \rangle_2 := \sum_{x \in G} u(x) v(x). \tag{5.13}$$

In the following theorem, we apply Lemma 2.5 to the construction laid out in [JP11]. This shows how one can recover the closability results described in [JP6, JP1, JP8, JP9] in a manner which is both quicker and more elegant. See also [JP10].

**Theorem 5.17** ([JP10, Thm. 5.22]). *Define* $K : \operatorname{span}\{\delta_x\}_{x \in G} \to \mathcal{H}_\mathcal{E}$ *by* $K\delta_x = \delta_x$ *and define* $\Delta$ *as in* (5.9). *Then* $(K, \Delta)$ *is a symmetric pair, and hence both operators are closable.*

**Theorem 5.18.** *The following are equivalent.*

(i) *For the Laplacian in* (5.10), *i.e., considered as an operator* $\Delta_\mathcal{E} : \mathcal{H}_\mathcal{E} \to \mathcal{H}_\mathcal{E}$ *with* $\operatorname{dom} \Delta_\mathcal{E} = \operatorname{span}\{v_x\}_{x \in G}$, $\Delta_\mathcal{E}$ *is* not *essentially self-adjoint. In other words, there is a nonzero* $\psi \in \mathcal{H}_\mathcal{E}$ *satisfying*
$$\psi \in \operatorname{dom} \Delta_\mathcal{E}^\star \quad \text{and} \quad \Delta_\mathcal{E}^\star \psi = -\psi. \tag{5.14}$$

(ii) *For the Laplacian in* (5.9), *i.e., considered as an operator* $\Delta_\mathcal{E} : \mathcal{H}_\mathcal{E} \to \ell^2(G)$ *with* $\operatorname{dom} \Delta = \operatorname{span}\{v_x\}_{x \in G}$, *there is a nonzero* $\psi \in \mathcal{H}_\mathcal{E}$ *satisfying*

$$\psi \in \operatorname{dom}(\Delta^\star K^\star) \qquad \text{and} \qquad \Delta^\star K^\star \psi = -\psi. \tag{5.15}$$

(iii) *The inclusions* $\overline{K} \subseteq \Delta^\star$ *and* $\overline{\Delta} \subseteq K^\star$ *are proper, i.e., there is a nonzero* $\psi \in \mathcal{H}_\mathcal{E}$ *satisfying*

$$\begin{bmatrix} \psi \\ K^\star \psi \end{bmatrix} \in \Gamma(K^\star) \ominus \Gamma(\Delta).$$

*Proof.* The equivalence (i) $\iff$ (ii) follows from the definitions of the respective adjoint operators. The equivalence (ii) $\iff$ (iii) is a direct application of Lemma 2.15 to the present situation. □

*Remark* 5.19. In light of Theorem 5.18 and [JP6, Prop. 4.9], which establishes that $\Delta_\mathcal{E}$ may not be essentially self-adjoint on $\mathcal{H}_\mathcal{E}$, it is clear that $(K, \Delta)$ is *not*, in general, a maximal symmetric pair.

**Lemma 5.20.** $\mathcal{H}arm \subseteq \operatorname{dom} \Delta_\mathcal{E}^\star$ *and* $\Delta_\mathcal{E}^\star h = 0$ *for all* $h \in \mathcal{H}arm$.

*Proof.* Let $u \in \operatorname{dom} \Delta_\mathcal{E}$ with $\Delta_\mathcal{E} u \in \mathcal{F}in$; for example, let $u = v_x$ for some $x \in G$. Then $\Delta_\mathcal{E} u \in \mathcal{F}in = \{\delta_x\}_{x \in G}^{\perp\perp}$ iff $\Delta_\mathcal{E} u \in \mathcal{H}arm^\perp$ by Theorem 5.11, so

$$\langle \Delta_\mathcal{E} u, h \rangle_{\mathcal{H}_\mathcal{E}} = 0. \tag{5.16}$$

It is clear this also shows $\Delta_\mathcal{E}^\star h = 0$. □

## References

[Aro] N. Aronszajn. Theory of reproducing kernels. *Trans. Amer. Math. Soc.* **68**(1950), 337–404. 14

[Bel] Denis Bell. The Malliavin calculus and hypoelliptic differential operators. *Infin. Dimens. Anal. Quantum Probab. Relat. Top.* **18**(2015), 1550001, 24. 2, 7

[BR] Ola Bratteli and Derek W. Robinson. *Operator algebras and quantum statistical mechanics. Vol. 1.* Springer-Verlag, New York, 1979. $C^*$- and $W^*$-algebras, algebras, symmetry groups, decomposition of states, Texts and Monographs in Physics. 2, 8, 11, 12

[DS] Nelson Dunford and Jacob T. Schwartz. *Linear operators. Part II.* Wiley Classics Library. John Wiley & Sons Inc., New York, 1988. 2, 3, 12

[FŌT] Masatoshi Fukushima, Yōichi Ōshima, and Masayoshi Takeda. *Dirichlet forms and symmetric Markov processes*, volume 19 of *de Gruyter Studies in Mathematics*. Walter de Gruyter & Co., Berlin, 1994. 14

[Gaw] Leszek Gawarecki. Transformations of index set for Skorokhod integral with respect to Gaussian processes. *J. Appl. Math. Stochastic Anal.* **12**(1999), 105–111. 2, 7

[Hid] Takeyuki Hida. *Brownian motion*, volume 11 of *Applications of Mathematics*. Springer-Verlag, New York, 1980. Translated from the Japanese by the author and T. P. Speed. 2, 7, 8

[HV] Cyril Houdayer and Stefaan Vaes. Type III factors with unique Cartan decomposition. *J. Math. Pures Appl. (9)* **100**(2013), 564–590. 2, 11, 12

[JP1] Palle E. T. Jorgensen and Erin P. J. Pearse. Operator theory and analysis of infinite resistance networks. (2009.), 1–247. arXiv:0806.3881. 13, 14, 15

[JP2] Palle E. T. Jorgensen and Erin P. J. Pearse. Unbounded containment in the energy space of a network and the Krein extension of the energy Laplacian. (2009.), 1–247. arXiv:1504.01332. 13, 15

[JP3] Palle E. T. Jorgensen and Erin P. J. Pearse. A Hilbert space approach to effective resistance metrics. *Complex Anal. Oper. Theory* **4**(2010), 975–1030. arXiv:0906.2535. 13, 15

[JP4] Palle E. T. Jorgensen and Erin P. J. Pearse. Resistance boundaries of infinite networks. In *Progress in Probability: Boundaries and Spectral Theory*, volume 64, pages 113–143. Birkhauser, 2010. arXiv:0909.1518. 13, 15

[JP5] Palle E. T. Jorgensen and Erin P. J. Pearse. Gel'fand triples and boundaries of infinite networks. *New York J. Math.* **17**(2011), 745–781. arXiv:0906.2745. 13, 15

[JP6] Palle E. T. Jorgensen and Erin P. J. Pearse. Spectral reciprocity and matrix representations of unbounded operators. *J. Funct. Anal.* **261**(2011), 749–776. arXiv:0911.0185. 13, 15, 16

[JP7] Palle E. T. Jorgensen and Erin P. J. Pearse. A discrete Gauss-Green identity for unbounded Laplace operators, and the transience of random walks. *Israel J. Math.* **196**(2013), 113–160. arXiv:0906.1586. 13, 14, 15

[JP8] Palle E. T. Jorgensen and Erin P. J. Pearse. Multiplication operators on the energy space. *J. Operator Theory* **69**(2013), 135–159. arXiv:1007.3516. 13, 15

[JP9] Palle E. T. Jorgensen and Erin P. J. Pearse. Spectral comparisons between networks with different conductance functions. *Journal of Operator Theory* **72**(2014), 71–86. arXiv:1107.2786. 13, 15

[JP10] Palle E. T. Jorgensen and Erin P. J. Pearse. Symmetric pairs and self-adjoint extensions of operators, with applications to energy networks. *To appear: Complex Anal. Oper. Theory* (2015), 1–11. arXiv:1512.03463. 1, 2, 3, 13, 15

[JP11] Palle E. T. Jorgensen and Erin P. J. Pearse. Unbounded containment in the energy space of a network and the Krein extension of the energy Laplacian. (2015). 17 pages, in review. arXiv:1504.01332. 15

[JPT] Palle E. T. Jorgensen, Erin P. J. Pearse, and Feng Tian. Duality for Unbounded Operators, and Applications. *In review* (2015), 1–14. arXiv:1509.08024. 2

[Kad] Richard V. Kadison. Dual cones and Tomita-Takesaki theory. In *Operator algebras and operator theory (Shanghai, 1997)*, volume 228 of *Contemp. Math.*, pages 151–178. Amer. Math. Soc., Providence, RI, 1998. 2, 11

[KR1] Richard V. Kadison and John R. Ringrose. *Fundamentals of the theory of operator algebras. Vol. II*, volume 100 of *Pure and Applied Mathematics*. Academic Press, Inc., Orlando, FL, 1986. Advanced theory. 2, 11

[KR2] Richard V. Kadison and John R. Ringrose. *Fundamentals of the theory of operator algebras. Vol. I*, volume 15 of *Graduate Studies in Mathematics*. American Mathematical Society, Providence, RI, 1997. Elementary theory, Reprint of the 1983 original. 3, 11

[Kat] Tosio Kato. *Perturbation theory for linear operators*. Classics in Mathematics. Springer-Verlag, Berlin, 1995. Reprint of the 1980 edition. 14

[KL1] Matthias Keller and Daniel Lenz. Dirichlet forms and stochastic completeness of graphs and subgraphs. *Preprint* (2009). arXiv:0904.2985. 13

[KL2] Matthias Keller and Daniel Lenz. Unbounded Laplacians on graphs: basic spectral properties and the heat equation. *Math. Model. Nat. Phenom.* **5**(2010), 198–224. 13

[LP] Russell Lyons and Yuval Peres. *Probability on Trees and Graphs*. Unpublished. 13, 15

[PS] K. R. Parthasarathy and K. Schmidt. *Positive definite kernels, continuous tensor products, and central limit theorems of probability theory*. Lecture Notes in Mathematics, Vol. 272. Springer-Verlag, Berlin, 1972. 2, 8

[Soa] Paolo M. Soardi. *Potential theory on infinite networks*, volume 1590 of *Lecture Notes in Mathematics*. Springer-Verlag, Berlin, 1994. 15

[Tak1] M. Takesaki. *Tomita's theory of modular Hilbert algebras and its applications*. Lecture Notes in Mathematics, Vol. 128. Springer-Verlag, Berlin-New York, 1970. 2, 11

[Tak2] M. Takesaki. *Theory of operator algebras. II*, volume 125 of *Encyclopaedia of Mathematical Sciences*. Springer-Verlag, Berlin, 2003. Operator Algebras and Non-commutative Geometry, 6. 2, 11

[vD] Alfons van Daele. The Tomita-Takesaki theory for von Neumann algebras with a separating and cyclic vector. In $C^*$-algebras and their applications to statistical mechanics and quantum field theory (Proc. Internat. School of Physics "Enrico Fermi", Course LX, Varenna, 1973), pages 19–28. North-Holland, Amsterdam, 1976. 2, 11

[Yam] Maretsugu Yamasaki. Discrete potentials on an infinite network. *Mem. Fac. Sci. Shimane Univ.* **13**(1979), 31–44. 15

University of Iowa, Iowa City, IA 52246-1419 USA,    palle-jorgensen@uiowa.edu

California Polytechnic University, San Luis Obispo, CA 93407-0403 USA,    epearse@calpoly.edu